\documentclass[10pt,a4paper]{article}
\usepackage{graphicx}
\usepackage{grffile}
\usepackage{longtable}
\usepackage{wrapfig}
\usepackage{rotating}
\usepackage[normalem]{ulem}
\usepackage{amsmath}
\usepackage{textcomp}
\usepackage{amssymb}
\usepackage{capt-of}
\usepackage[colorlinks=true]{hyperref}
\usepackage{minted}
\usepackage[margin=1.5cm]{geometry}
\usepackage{amsmath,amssymb}
\usepackage{graphicx}
\usepackage{booktabs}
\usepackage{longtable}
\usepackage{tabularx}
\usepackage{array}
\usepackage{enumitem}
\usepackage[round,authoryear]{natbib}
\usepackage{hyperref}
\hypersetup{colorlinks=true,linkcolor=black,citecolor=black,urlcolor=blue}
\setlist{nosep,leftmargin=*}

\newcolumntype{Y}{>{\raggedright\arraybackslash}X}
\newenvironment{dialogue}{%
\begin{list}{}{%
\setlength{\leftmargin}{2.2cm}%
\setlength{\labelwidth}{1.9cm}%
\setlength{\labelsep}{0.2cm}%
\setlength{\itemsep}{0.45\baselineskip}%
\setlength{\parsep}{0pt}%
\setlength{\topsep}{0.35\baselineskip}%
}%
}{\end{list}}
\newcommand{\speaker}[1]{\item[\textbf{#1}]}
\author{Siniša Miličić (ORCID: \href{https://orcid.org/0000-0002-1825-0097}{0000-0002-5659-5598})}
\date{}
\title{Open Preparation, Human Explanation, and Instructor Synthesis\\\medskip
\large A Human-Scale Methodology for AI-Rich Higher Education}
\hypersetup{
 pdfauthor={Siniša Miličić (ORCID: \href{https://orcid.org/0000-0002-1825-0097}{0000-0002-5659-5598})},
 pdftitle={Open Preparation, Human Explanation, and Instructor Synthesis},
 pdfkeywords={},
 pdfsubject={},
 pdfcreator={Emacs 30.2 (Org mode 9.7.30)}, 
 pdflang={English}}
\begin{document}

\maketitle
\section{Abstract}
\label{sec:org81a9741}

In AI-rich higher education, polished written mathematics has become easier to
produce than trustworthy evidence of understanding. This article develops a
human-scale methodology for service mathematics, with informatics as its main running
case. Its central move is not prohibition of tools but relocation of evidential
trust. Students prepare openly, often with digital assistance, but grade-relevant
evidence shifts toward live explanation, contingent questioning, and cumulative
observation against course outcomes. The design is guided by Realistic Mathematics
Education, question-first task construction, short human-scale mathematical tasks,
and instructor synthesis after student attempt. It contributes a weekly operational
cycle, a realism framework distinguishing professional, disciplinary, and
experiential realism, a middle-out white-box / black-box stance on tools, a bounded
role for retrieval-grounded AI assistants for students and teachers, and a cumulative
oral-evidence model for small and medium cohorts. The paper also provides concrete
implementation artifacts: process figures, an ecology of problem types, time-budget
estimates, an evidence hierarchy, and a five-grade oral rubric. This is a methodology
paper rather than an effectiveness study. Its claim is that the proposed design is
pedagogically coherent, operationally plausible for human-scale teaching settings,
and responsive to current concerns about AI, oral evidencing, and active learning in
undergraduate mathematics education.

\medskip\noindent\textbf{Keywords:} service mathematics; Realistic Mathematics Education; oral assessment; generative AI; flipped classroom; undergraduate mathematics education
\section{Introduction}
\label{sec:org3bf078c}

Three pressures now meet in the same classroom.

\begin{enumerate}
\item University mathematics still needs coherence, not only coverage.
\item Service mathematics still needs students to explain, justify, and connect formal
ideas.
\item Generative AI has made polished written output cheap in comparison with
trustworthy evidence of understanding.
\end{enumerate}

This paper answers those pressures with a methodology rather than a ban. Its claim is
that service mathematics should be organized as a question-first, RME-guided, orally
validated learning design: students meet short coherent problems before class,
prepare openly, explain publicly, are questioned contingently, and only then receive
fuller instructor synthesis. Written artifacts remain useful, but they lose their
privileged status as self-authenticating evidence.

The problem is not AI alone. The deeper educational anti-standard is the fantasy that
difficult intellectual work can be delegated to a machine, skipped as a human
journey, and then somehow recovered at the destination. Recent discussion of AI and
mathematics warns against being dropped at the destination while missing the trail,
the intermediate structure, and the human gains of the journey, and separately
criticizes push-of-a-button workflows that reduce mathematical activity to button
pressing rather than structured human--AI conversation \citep{Wong2026TaoAI}. The
present paper takes that warning seriously, but translates it into undergraduate
service mathematics. Students may prepare with tools. They still have to explain as
humans.

The method is not anti-technology. It is humanly mathematical. Mathematics here is
treated as a deeply human activity whose educational role is to enable disciplined
reasoning about structures, representations, systems, and professional
problems. Informatics serves as the main running case because it makes the present
evidential problem especially sharp: it is mathematically adjacent, but not
mathematics; it is structurally rich, yet often weakly constrained by physical
immediacy. At the same time, the method is not limited to informatics. It has natural
relevance wherever mathematics is taught in support of adjacent professional fields.

\begin{figure*}
\centering
\includegraphics[width=0.98\textwidth]{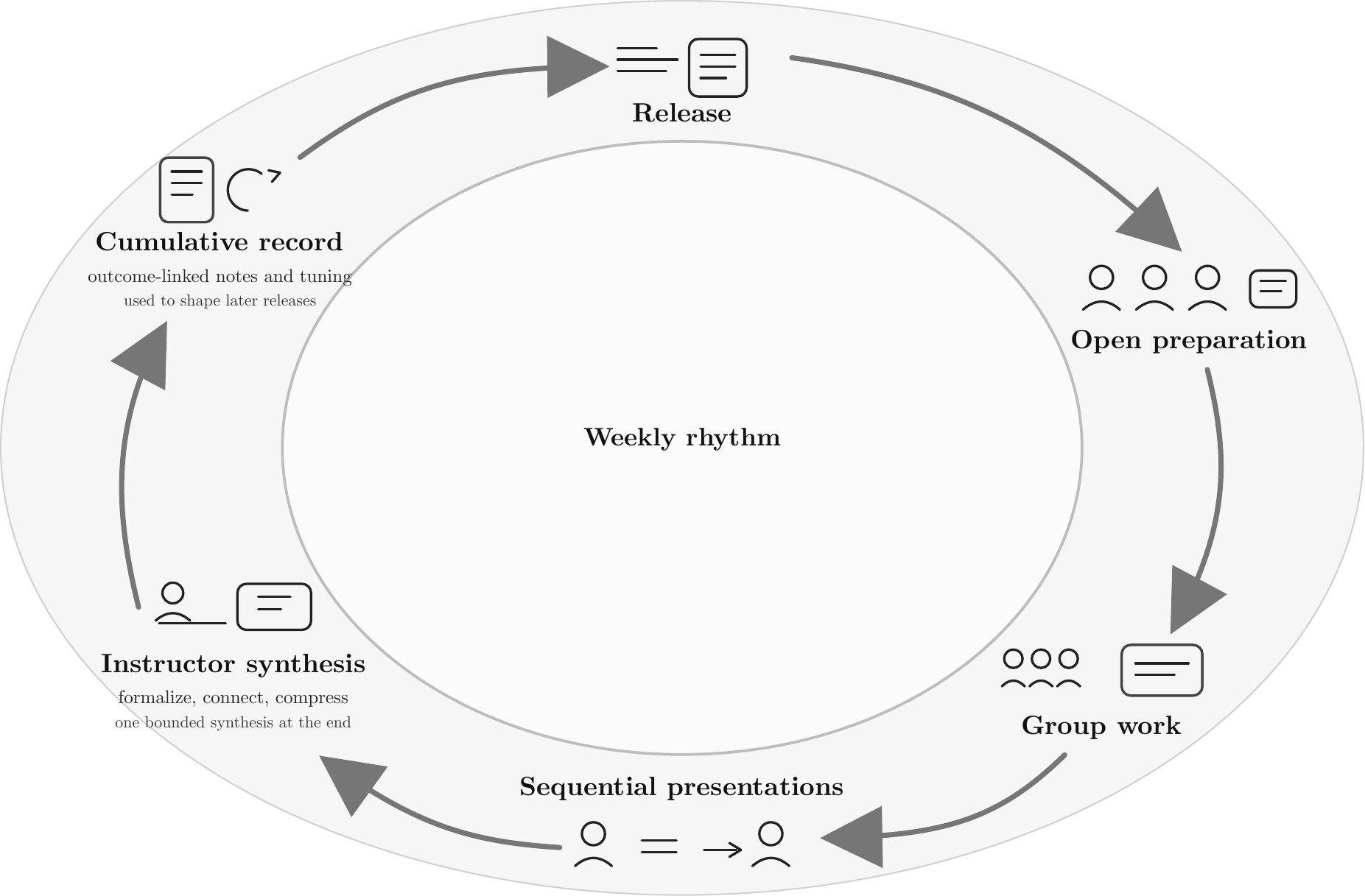}
\caption{\label{fig:orgdf1d789}The figure shows the recurring weekly rhythm of the proposed question-first, orally validated design.}
\end{figure*}

Figure \ref{fig:orgdf1d789} introduces the operational contrast that motivates the rest
of the paper. The first gives a compact sketch of the conventional
lecture--homework--exam process. The second gives the full weekly process proposed
here. Together they make the shift visible: first contact no longer resides in
lecture, synthesis no longer precedes student attempt, and validity no longer resides
mainly in the final written artifact.

A deliberate feature of the manuscript is that it partly applies its own method to
itself. It begins with the practical problem, moves quickly to one realistic
situation, and only then turns to broader pedagogical analysis. In that limited
sense, the article is not only about a question-first method; it is partially
structured by it.

The argument proceeds in six steps. First, a running example is introduced. Second,
the pedagogical commitments are stated in relation to RME, Kurnik, IBL / PBL, and a
middle-out white-box / black-box stance. Third, the weekly operational model and the
two Oracles are specified. Fourth, realism, an ecology of problem types, and oral
ownership are developed. Fifth, evidence, time, and grading are formalized. Sixth,
the paper closes with a short discussion of limits, neighboring design languages, and
future research.
\section{A running example: the grayscale operator}
\label{sec:org902cbbf}

Because the method is question-first, the paper should not begin only with
abstractions. It therefore opens its pedagogical argument through a single running
example from geometry and linear algebra.

A student group is asked to analyze a grayscale transformation on RGB pixels. A pixel
in color is represented as a vector \(c=(R,G,B)\). The grayscale procedure sends each
pixel to a gray value \(s=0.3R+0.6G+0.1B\), repeated across all three
coordinates. The immediate student task is not ``write down the matrix because linear
algebra says so,'' but something more local and intelligible: determine what happens
to pure red, pure green, and pure blue; compare the outputs; and decide what kind of
transformation this seems to be.

That first phase is already mathematical. Students begin from a familiar
computational act, identify stable images of the basis colors, and notice that all
outputs land on the grayscale line in \(\mathbb{R}^{3}\). Only later does the
synthesis formalize the process as a linear operator, represent it by a matrix, and
connect repeated application with idempotence and projection-like behavior.

This example is useful throughout the paper for three reasons. It is experientially
accessible, professionally plausible within computing, and mathematically rich enough
to illustrate basis images, operator representation, structural closure, and the
middle-out handling of white and black boxes. It therefore serves as the paper's main
running case.

The example also helps make an accessibility point that matters more broadly for the
methodology. RGB technology is now common cultural knowledge even for many readers
who are color-blind. The pedagogical point here is not color sensation as such, but
the structure of the operator and the movement from situation to
formalization. Color-dependent examples should nevertheless be chosen and presented
so that their essential structure remains legible to color-blind students and to
readers encountering figures in grayscale print or on grayscale e-ink screens.
\section{Pedagogical positioning}
\label{sec:org99d686f}

The method rests on five commitments.

\begin{enumerate}
\item Mathematics should begin as guided activity rather than finished reception.
\item Coherence should be authored question-first rather than topic-first.
\item Realistic contexts should be chosen for students, not only for teachers.
\item Tool use should be handled through a middle-out white-box / black-box pedagogy.
\item Oral explanation should be treated as constitutive of learning, not only as evidence of learning.
\end{enumerate}

The first commitment is Freudenthalian. Realistic Mathematics Education treats
mathematics as a human activity and asks that learners encounter mathematically
meaningful situations from which concepts, structures, and methods may emerge through
guided reinvention \citep{Freudenthal1991Revisiting,HeuvelPanhuizenDrijvers2014RME,GravemeijerDoorman1999}. The present method is not a doctrinally pure RME program,
but it is strongly shaped by RME's view of mathematics as a human activity, by guided
reinvention, and by the demand that contexts be realistic for learners rather than
only for teachers.

The second commitment is architectural. The weekly design starts from questions. This
paper does not carry the full topical-unit theory as its main burden, but it is still
informed by it: a good weekly problem should be small enough to be thinkable,
teachable, orally defensible, and later recombinable into larger structures. In
practice, question-first design serves as a discipline of coherence: what exactly is
this week trying to make answerable?

\begin{figure*}
\centering
\includegraphics[width=0.90\textwidth]{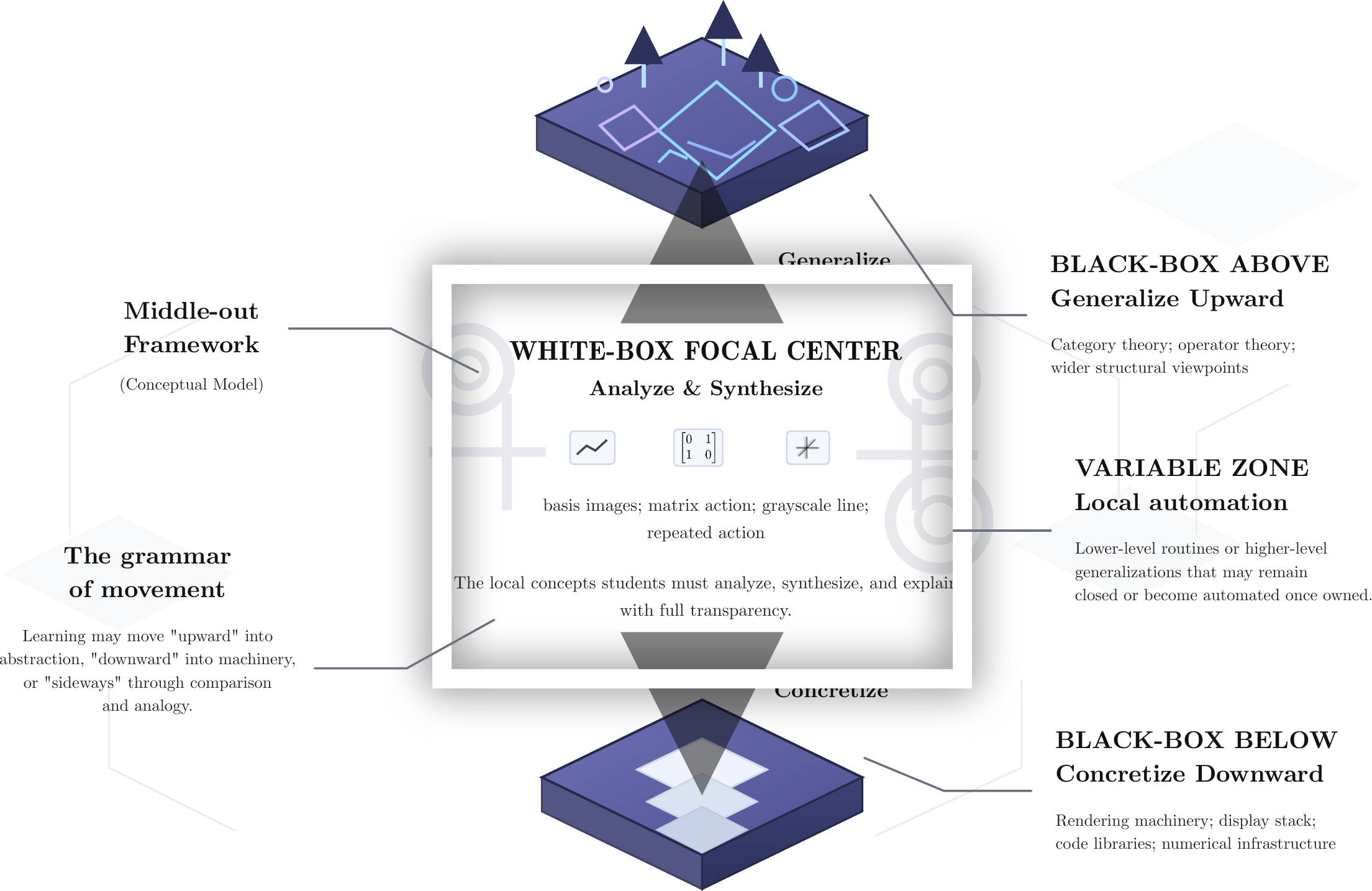}
\caption{\label{fig:orgdabac4e}Conceptual middle-out structure. The focal concept sits among deeper machinery, wider abstraction, and locally variable structure.}
\end{figure*}

The third commitment requires care. A context can be highly realistic from a
teacher's institutional or professional perspective while remaining experientially
empty to students. A familiar higher-education computing cliche is course timetabling
and room allocation: it is administratively real, algorithmically rich, and
intellectually defensible, but for many first-year students it is experienced as the
university's scheduling problem rather than as their own mathematical situation
\citep{Johnson1993Timetabling,MullerBartak2002InteractiveTimetabling}. More
generally, service mathematics can drift into contexts that are realistic to the
institution, to academia, or to a profession, without being experientially real to
the student. RME is a useful corrective because it asks whose reality is being
invoked.

\begin{table*}[htbp]
\caption{\label{tab:org6e13264}Middle-out adaptation of the white-box / black-box distinction.}
\centering
\begin{tabularx}{\textwidth}{lYY}
Layer & Typical movement & Grayscale example\\
\hline
White-box & Analyze; synthesize & Basis images; matrix action; grayscale line; repeat action\\
Variable & Open or close locally & Routine matrix work; fast basis-to-matrix translation\\
Black-box below & Concretize downward if needed & Renderer; display stack; viewer internals; color pipeline; numeric libraries\\
Black-box above & Generalize upward if needed & Jordan form; category theory; operator theory\\
\end{tabularx}
\end{table*}

The fourth commitment needs a more careful statement than the familiar ``learn the
basics, then use tools.'' That linear slogan is pedagogically weak. It suggests that
mathematical movement in a course has one privileged direction, either upward from
foundations or downward from abstractions. In practice, neither description is
adequate. A teacher who insists that a course is purely ``bottom-up'' often
underestimates how much wider structure already governs the local task; a teacher who
insists that it is purely ``top-down'' often underestimates how much conceptual work
students must still build underneath the current object. In an actual course, there
is always something ``under'' the present concept, something ``over'' it, and
something conceptually adjacent to it. For that reason, one always begins from a
middle.

I therefore adopt an explicitly \emph{middle-out} stance. The point is not compromise
between two directions, but recognition that mathematically meaningful teaching
proceeds from a local conceptual middle and then moves in several possible directions
as needed. Those directions include analysis and synthesis, abstraction and
concretization, generalization and specialization, induction and deduction, and also
comparison and analogy. In Kurnik's methodological language, these belong to the
basic scientific methods relevant to mathematics teaching; in the present framework,
they function as a practical grammar of movement between concepts rather than as a
fixed instructional sequence \citep{Kurnik2008Znanstvenost,Kurnik2009Problemska}.

On this view, ``white-box'' and ``black-box'' are not simply endpoints on a one-way
developmental ladder. They mark temporary didactical decisions about which
surrounding structures must be opened, which may remain closed, and which may shift
status as students gain ownership. Some movements go downward into underlying
mechanism, some upward into wider abstraction, some sideways through comparison and
analogy, and some inward through analysis and later outward through synthesis around
the focal object. The grayscale example shows this clearly: basis images, matrix
action, and repeated application should be opened locally; rendering machinery may
remain closed below; operator-theoretic and categorical generalizations may remain
closed above; and routine calculations may move between the two depending on observed
mastery. Category theory is one plausible ``outer'' direction here, not because the
local task is secretly about category theory, but because a local operator-centered
discussion may eventually point toward functors and more general structural
viewpoints \citep{MacLane1998Categories}. At the same time, categorical ideas may
also appear ``below'' the present teaching target inside programming-language
semantics, compiler design, or software infrastructure, even when the local course
problem does not open them explicitly \citep{Pierce1991BasicCategoryTheory}. The
stance is therefore middle-out not because it avoids direction, but because it
recognizes a richer directional grammar than the simple opposition of ``bottom-up''
and ``top-down.''

Table \ref{tab:org6e13264} and Figure \ref{fig:orgdabac4e} summarize this stance in
complementary ways: the table lists the didactical status of surrounding structures,
while the figure shows the focal concept as situated within a wider grammar of
possible conceptual movement.

The fifth commitment concerns evidencing, but also learning itself. Oral explanation
matters here not because it is old-fashioned or theatrical, but because it remains
one of the clearest places where mathematical understanding becomes visible under
follow-up questioning. More than that, presenting and answering are themselves
constitutive of learning: students do not merely display understanding there, they
further organize it there. That constitutive role and that evidential role should be
held together from the beginning. The method aligns them rather than treating them as
separate concerns. This brings it into contact with flipped learning,
peer-instruction-like cycles, constructive alignment, and cumulative assessment
\citep{Strelan2020FlippedMeta,VanAlten2019FlippedMeta,LoHewChen2017MathFlipped,Biggs1996Constructive,VanderVleutenEtAl2012Programmatic,BaartmanVanSchiltMolVanderVleuten2022}. It also stands in continuity with a Croatian
methodological lineage in which mathematics teaching is organized around problem
situations, heuristic guidance, and the scientific organization of student
inquiry\footnote{Zdravko Kurnik's writings on \emph{problemska nastava} and \emph{znanstvenost u
nastavi matematike} are especially relevant here. Their common thesis is that
mathematical understanding should emerge through structured work on problems, guided
questioning, and progressive generalization rather than passive reception of finished
results. The present framework stands in continuity with that problem-centered
lineage, while extending it toward explicit course architecture, cumulative oral
evidencing, and AI-aware validity \citep{Kurnik2009Problemska,Kurnik2008Znanstvenost}.}.

This also helps distinguish the present method from standard IBL and PBL. It is
positively related to both, but it is more tightly orchestrated. Compared with
generic IBL, it builds inquiry into a weekly cycle that explicitly couples
problem-first work with oral validation and instructor synthesis. Compared with
classic PBL, it does not revolve around large multi-session scenarios; it relies on
short human-scale mathematical tasks that can adapt to attendance variation and to
observed levels of understanding \citep{HmeloSilver2004PBL,LaursenRasmussen2019IBL}.
\section{Realism in this method}
\label{sec:org2319b61}

\begin{figure*}
\centering
\includegraphics[width=0.98\textwidth]{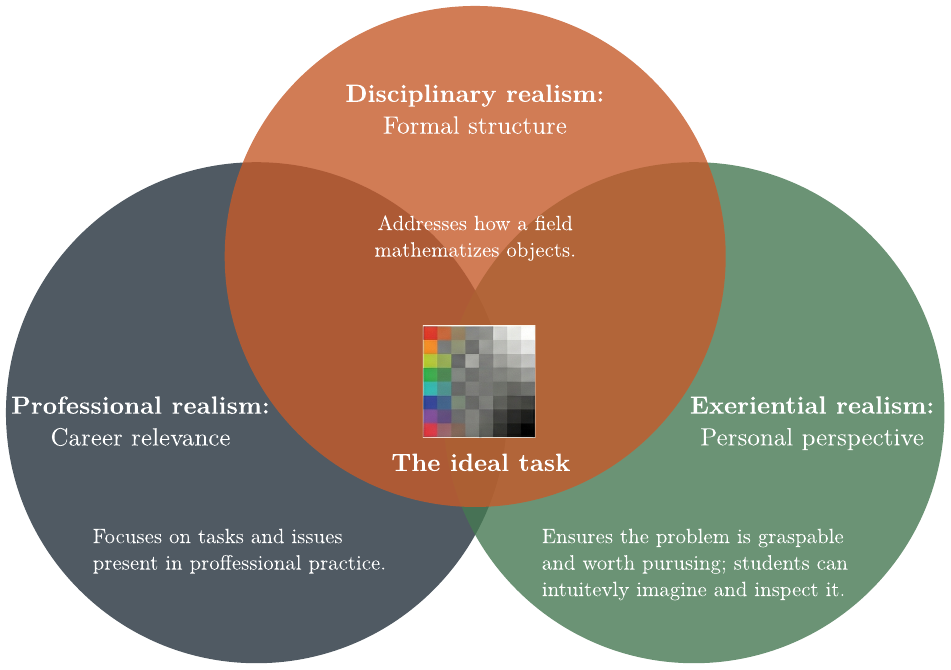}
\caption{\label{fig:orgc05c0e1}Professional, disciplinary and experiental realisms.p}
\end{figure*}

Three kinds of realism matter here, represented in Figure \ref{fig:orgc05c0e1}.

\begin{enumerate}
\item \emph{Professional realism}: the task is realistic relative to actual or near-future
professional practice that uses structure, representation, constraint, or formal
procedure, whether or not practitioners name that work as mathematics.
\item \emph{Disciplinary realism}: the task is realistic relative to how a discipline
structures, idealizes, and mathematizes its own objects.
\item \emph{Experiential realism}: the task is graspable, imaginable, and worth acting on
from the learner's point of view.
\end{enumerate}

These categories should not be collapsed. Professional realism concerns the world of
career practice: software systems, protocols, files, interfaces, signals, data, and
computation. Disciplinary realism concerns the inner formal framing of a field: what
counts as an object, what structural features matter, what formalization is natural,
and how explanations are organized. Experiential realism concerns the learner's
actual entry point: whether the task is tangible enough, interpretable enough, and
worth pursuing from the student's perspective.

A short clarification helps here. By \emph{mathematizing}, I mean the act of recasting a
practical or disciplinary situation in terms of objects, relations, constraints,
transformations, or invariants that can be reasoned about systematically. In
undergraduate informatics, this may happen even where participants insist they are
not ``doing mathematics.'' Consider protocol formatting: a message specification may
initially appear as a practical implementation detail, but once one begins to reason
about field order, admissible values, encoding constraints, compositional structure,
or lossless transformations between representations, the activity has already moved
into mathematizing. The point is not that every such activity should be renamed
``mathematics,'' but that disciplinary realism often depends on these structuring
moves whether or not they are explicitly recognized as such.

The grayscale task is again helpful. It is professionally real because color handling
and image processing belong recognizably to computing practice. It is disciplinarily
real because the relevant structure is genuinely linear-algebraic: basis images,
target subspaces, and repeated action. It is experientially real because nearly every
student can imagine what it means to turn color into gray and can inspect what
happens to primary colors.

This is also where the present method earns the RME label more concretely. The
question-first phase is not merely pre-lecture homework. It is meant to begin a local
process of mathematization. Students first organize a situation, identify a stable
pattern, and only then receive fuller formal notation and synthesis. In the grayscale
task, for example, students can begin from the familiar act of turning a colored
pixel into gray, inspect what happens to pure red, green, and blue, notice that all
outputs lie on the grayscale line, and only then receive the formal operator language
in synthesis. The task is realistic not only because it has a computational context,
but because it lets students move from situation to structure to formalization.

A good task therefore does not merely ``apply mathematics to something real.'' It
chooses a context whose reality works at the right level for this student group, at
this moment, in this course, and then supports mathematization rather than merely
attaching a formula to an application story.

\section{The weekly operational model}
\label{sec:orgdcca1dd}

The weekly cycle is the practical core of the method. Figure \ref{fig:orgdf1d789} gives
the full proposed process, while Table \ref{tab:org14e173b} states the same cycle in
compact teacher-usable terms. The figure carries some ontology of scale that is worth
making explicit in prose as well. Release, synthesis, and cumulative recording occur
once per weekly cycle. Preparation occurs per student and in parallel. Group work
occurs in a smaller number of groups and in parallel. Presentations and follow-up
occur sequentially. The cycle then feeds the next release. Stating this directly is
useful because the human-time arithmetic of the method depends on it.

One presented problem is budgeted at about five to seven minutes plus transition. A
single weekly cycle then includes one final synthesis segment, typically amounting to
about three minutes per presented task in aggregate rather than three minutes
inserted after each presentation. That distinction matters for the arithmetic.

The in-class core should be read very plainly. Students are assigned to groups truly
randomly. Each group receives one problem. Groups get about 20 minutes of breakout
preparation. One person from each group presents. Group formation is procedurally
fair by mechanical randomness. Presenter selection is humanely fair: it is
randomized, but favors students who have been seen less often presenting and remains
subject to humane teacher adjustment for directly visible circumstances such as
panic, temporary setback, or other immediately recognizable situations. Acquaintance
and audience familiarity come from the recurring cohort, not from fixed
groups. Synthesis is the final integrative part.

That final humane clause matters. Kindness and humanity come first. The method is
structured enough to resist caprice, but not so rigid that it overrides obvious
pedagogical judgment.

\begin{table*}[htbp]
\caption{\label{tab:org14e173b}Operational description of the weekly cycle.}
\centering
\begin{tabularx}{\textwidth}{lYYY}
Phase & Action & Role & Emphasis\\
\hline
Release & Problem bundle & First encounter & RME prompt; question-first frame\\
Preparation & Notes; readings; code; AI & First sense-making & Inquiry; problem engagement\\
Group work & Random groups; one problem; 20 min & Clarify; compare & Peer instruction; mathematization\\
Presentation & One speaker / group & Main visible evidence & Constitutive learning via assessment\\
Follow-up & Teacher + class probe & Ownership check & Repair; adaptation; structural check\\
Synthesis & Integrate; correct; formalize & Closure & Feedback; compression; formalization\\
Record & Log by outcome & Longitudinal judgment & Validity + formative continuity\\
Supplement & Handout / short video & Post-class support & Not a replacement for live work\\
\end{tabularx}
\end{table*}
\subsection{Week 1 as bootstrapping}
\label{sec:org24e3237}

Week 1 is structurally different. Students cannot yet pre-prepare in the normal way,
so the course begins in \emph{bootstrapping mode}. This week has three purposes.

\begin{enumerate}
\item lowered-stakes rehearsal of the cycle,
\item explicit norm-setting for the semester,
\item early modeling of what counts as good oral mathematical explanation.
\end{enumerate}

Students still present in week 1, but the problems are intentionally smaller, the
stakes lower, and the teacher more interventionist. The point is not early
sorting. It is to establish the rhythm in a humane way before the rhythm becomes
consequential.
\subsection{Synthesis duration and function}
\label{sec:orgefd7ddb}

The synthesis phase is not a decorative wrap-up. In a standard weekly meeting it is a
bounded final segment. A useful rough estimate is about three minutes of synthesis
per presented problem, but that synthesis is grouped at the end rather than
distributed after each student presentation. The phase remains variable and
topic-sensitive, but its function is stable.

\begin{enumerate}
\item It formalizes and compresses what students have already tried to articulate.
\item It models disciplined mathematical explanation for students and, where relevant,
for less experienced teachers.
\end{enumerate}

In the grayscale case, after a student has presented the images of the basis vectors
and described the operator once in computational terms, the synthesis can step back
and say: what we have really built is not just a recipe for turning colored pixels
gray, but a linear map with a specific target structure. Once a pixel is already on
the grayscale line, applying the operator again changes nothing further; that is why
the transform behaves like an idempotent projection onto that line, and why basis
images tell us so much about the whole map.
\subsection{Anxiety, language, and the affective transition}
\label{sec:org3e65d70}

Because the method shifts grade-relevant evidence toward oral work, its affective
design matters. The general principle is that each particular oral session is
low-pressure and partial-credit in nature. A student is not expected to carry the
whole course or the whole topic in one spotlight performance. The cumulative
structure matters precisely because each local oral moment is limited, sampled, and
repeatable.

This helps with communication anxiety, with students for whom the language of
instruction is not their strongest language, and with students who need time to grow
into public mathematical speech. Frequency can reduce anxiety because public
explanation becomes normal rather than exceptional. The recurring local cohort
stabilizes the audience as well. The bootstrapping transition from week 1 to week 2
should therefore be read as a shift from rehearsal to accountability, but not from
kindness to harshness. Where needed, humane adjustments remain legitimate during the
semester, and unresolved gaps or borderline cases may be repaired through a
supplementary oral opportunity late in the semester rather than through routine
deferment.
\section{The Learning and Governance Oracles}
\label{sec:org22ebf48}

The method includes two grounded AI supports.

\begin{enumerate}
\item the \emph{Learning Oracle} for students,
\item the \emph{Governance Oracle} for teachers.
\end{enumerate}

These are not generic chatbots and should not be read as deus ex machina
solutions. They are bounded retrieval-augmented generation systems over course
materials and related corpora. In less technical language, they are RAG-enhanced LLM
wrappers over a limited document universe. On the student side, a NotebookLM-like
grounded system is one workable example; Google's own education materials describe
NotebookLM as grounded in the sources provided by the user and able to return in-line
citations, which makes it a reasonable model for the kind of bounded assistant
imagined here \citep{GoogleNotebookLM2026}. On the teacher side, the point is not
automation for its own sake, but governance of the weekly cycle. Educational RAG
literature already notes both the appeal and the limits of such systems: grounding
can reduce error and improve domain relevance, but it does not eliminate
hallucination or flawed reasoning \citep{LiEtAl2025RAGEducation}.

Table \ref{tab:org0dede6c} states the operational protocol.

\begin{table*}[htbp]
\caption{\label{tab:org0dede6c}Oracle protocol.}
\centering
\begin{tabularx}{\textwidth}{lYY}
Aspect & Learning Oracle & Governance Oracle\\
\hline
Corpus & Notes; tasks; worked examples; outcomes; readings; FAQ / misconceptions & Notes; problem bank; pedagogy sources; misconception bank; prior tasks\\
Main use & Prepare; reformulate; explain; navigate & Generate; vary; tune; aggregate\\
Hard limit & Provisional; must be checked & No student-facing live evaluation\\
\end{tabularx}
\end{table*}

The Learning Oracle is useful precisely because it is grounded, but it is not
authoritative. Grounded systems can still distort, oversimplify, or mislead. For that
reason, Oracle criticism is part of student responsibility: Oracle output is always
provisional and should be human-verified against course materials where
appropriate. The oral phase is the final filter for Oracle-borne misunderstanding. A
student who prepares with the Learning Oracle is still responsible for criticism,
verification, and oral ownership of the resulting explanation. In that sense,
Oracle-derived mistakes can be pedagogically productive if they are exposed,
criticized, and repaired in follow-up rather than treated deferentially. Teachers
should also watch for \emph{Oracle-borne misconceptions} as a distinct error class: they
are often polished, plausible, and locally coherent, which makes live follow-up
especially important.

The Governance Oracle supports task generation, adaptive variation, misconception
tracking, and weekly diagnosis from the teacher side. Its strongest function is
governance: keeping tasks aligned, varied, level-appropriate, and responsive to the
actual student cohort. In smaller groups, this creates a real agile advantage. The
course can be tuned week by week without losing coherence. A typical use case would
be: after week 4, the teacher logs a recurrent misconception in the cohort's
explanations of entropy or linear dependence, and the Governance Oracle aggregates
those diagnosed patterns and proposes three short task variants for the next release
that keep the same outcome target but alter context, representation, and difficulty.

A hard privacy rule is needed. The Governance Oracle belongs on the teacher side
only. It should work on de-identified or cohort-level patterns, not on named-student
live assessment support. Misconception tracking belongs primarily at cohort and topic
level; it should not become a mechanism for pigeonholing individual students in the
moment.
\section{Ecology of problem types}
\label{sec:org17481a0}

Table \ref{tab:org431881b} sketches the intended ecology of problem types. The point
is not one flagship gimmick, but a deliberately mixed set of small problem types that
are rich enough for follow-up, realistic in more than one sense, and suitable for
question-first work. The running grayscale case appears here again, but now among
peers rather than alone. The additional functional-programming example is included
because it reveals especially clearly how the same object may be treated locally as a
programming abstraction, upward as a category-theoretic direction, and downward as
part of semantic or infrastructural machinery not opened in the present task. Similar
logic also extends beyond the present informatics center of gravity, for example to
weather and sensor structure in maritime studies, geometric modeling in mechanical
engineering, or data-science settings shared across many current programs.

The grayscale operator is especially useful in this ecology because it lets one see,
in a single small task, the whole movement from practical recognition, to
mathematization, to contingent structural questioning.

\begin{table*}[htbp]
\caption{\label{tab:org431881b}Ecology of problem types.}
\centering
\begin{tabular}{>{\raggedright\arraybackslash}p{0.12\textwidth} >{\raggedright\arraybackslash}p{0.14\textwidth} >{\raggedright\arraybackslash}p{0.17\textwidth} >{\raggedright\arraybackslash}p{0.17\textwidth} >{\raggedright\arraybackslash}p{0.17\textwidth} >{\raggedright\arraybackslash}p{0.15\textwidth}}
Family & Topic & Professional realism & Disciplinary realism & Experiential realism & Guiding question\\
\hline
Information / coding & Kolmogorov + compression & Files; compression & Structure and randomness & Run code; inspect files & Why does one string compress but another resist?\\
Discrete metrics & Hamming vs Levenshtein & Transmission & Metric assumptions & Compare short words & What disturbance is this metric built to detect?\\
Linear algebra & Grayscale operator & Image processing & Basis images; subspace; repeat action & Familiar color cases & What kind of transform is this, and what does repeat action do?\\
Vector geometry & 3D triangle area & Graphics; simulation & Cross product as structure & Sketchable triangle & How does one vector encode area?\\
Graph theory & Connectivity failure & Networks & Reachability and proof & Small drawable graphs & What exactly breaks when connectivity fails?\\
Functional programming & Functor in FP & Typed data transformation & Structure-preserving mapping; composition & Concrete type constructors & What is preserved when a mapping operation becomes a functor?\\
\end{tabular}
\end{table*}

These topics are not unified by subject matter alone. They are unified by their
ability to support mathematization, oral explanation, contingent follow-up, and later
synthesis.
\section{A hypothetical oral dialogue: the grayscale operator, mathematization, and ownership}
\label{sec:org42e2884}

Listing \ref{lst:grayscale-dialogue} illustrates the difference between executing a
familiar computational routine and owning the mathematical structure that organizes
it. The setting is intentionally realistic: the student is a first-year
non-mathematics student who recognizes grayscale conversion from computing practice,
but has not yet fully mathematized it as a linear operator. The dialogue therefore
moves from a familiar situation, to structural interpretation, to contingent
questioning about what the operator preserves and what repeated application reveals.

\begin{listing*}[h!]
\centering
\begin{minipage}{0.99\textwidth}
\caption{A hypothetical follow-up dialogue on the grayscale operator. The aim is to expose the difference between local procedural description and owned structural explanation across computational, algebraic, and geometric views.}
\label{lst:grayscale-dialogue}
\begin{dialogue}

\speaker{Instructor} Tell me first, in ordinary computing language, what this grayscale procedure does.

\speaker{Student} It takes a color pixel and turns it into a shade of gray.

\speaker{Instructor} Good. Now say that more carefully. What is the input and what is the output?

\speaker{Student} The input is an RGB vector, so something like $(R,G,B)$. The output is again a vector, but all three entries become equal.

\speaker{Instructor} Equal to what?

\speaker{Student} To a weighted sum, something like $0.3R + 0.6G + 0.1B$.

\speaker{Instructor} So is this just an image-processing recipe, or is it also a mathematical object?

\speaker{Student} It is also a linear map from $\mathbb{R}^3$ to $\mathbb{R}^3$.

\speaker{Instructor} Why are you entitled to call it a linear map?

\speaker{Student} Because each output coordinate depends linearly on $R$, $G$, and $B$, and the same formula works uniformly for every input pixel.

\speaker{Instructor} Good. Then why do we bother looking at pure red, pure green, and pure blue?

\speaker{Student} Because they are the basis vectors. If we know where those go, we know the whole map.

\speaker{Instructor} That sounds like a slogan. Why is it true?

\speaker{Student} Every color vector can be written as a linear combination of those basis vectors, and a linear map is determined by what it does to linear combinations. So once we know the images of the basis vectors, we can rebuild the image of any pixel.

\speaker{Instructor} Excellent. Now say something geometric. Where do all outputs live?

\speaker{Student} On the grayscale line, where all three coordinates are equal.

\speaker{Instructor} And what does that tell you about the map?

\speaker{Student} It collapses the whole RGB space onto that line.

\speaker{Instructor} Suppose I apply the same grayscale operator again to a pixel that is already gray. What happens?

\speaker{Student} Nothing changes.

\speaker{Instructor} Why not?

\speaker{Student} Because once the vector is already on the grayscale line, the weighted sum gives back the same gray value. So applying the operator again keeps it in the same place.

\speaker{Instructor} So what kind of behavior is that?

\speaker{Student} It behaves idempotently, and in this setting it can be understood as a projection onto the grayscale line.

\speaker{Instructor} Good. What changed between your first answer and your last answer?

\speaker{Student} At first I described what the software does. Now I can explain the mathematical structure: the input space, the basis images, the target line, and why repeating the operation does not change the result.

\end{dialogue}
\end{minipage}
\end{listing*}

The movement is from practical recognition, to mathematization, to contingent
structural questioning. That is the kind of ownership the method is trying to make
visible: the student does not merely recognize the task from experience, but can
reorganize it as a mathematical object and defend that reorganization under
follow-up.
\section{Human-time budget and cohort size}
\label{sec:org478f8c5}

Any serious methodology paper must show the arithmetic of the classroom. The present
method is therefore stated natively for small and medium cohorts, not for industrial
scale. Human explanation and human judgment do not scale by magic.

For a two-session weekly cycle of \(2 \times 45\) minutes, the gross budget is 90
minutes. About 20 minutes are used for group work, leaving about 70 minutes for
student presentations and the final synthesis. If one task consumes about 10 minutes
of that remaining budget in gross terms---roughly 5 to 7 minutes of presentation plus
transition and about 3 minutes of synthesis allocation---the weekly cycle supports
about 6 to 8 student presentations together with about 20 to 25 minutes of final
synthesis. For a four-session weekly cycle of \(4 \times 45\) minutes, the same logic
leaves about 160 minutes after group work, allowing up to about 16 presentations
before student attention begins to degrade.

\begin{figure*}
\centering
\includegraphics[width=0.98\textwidth]{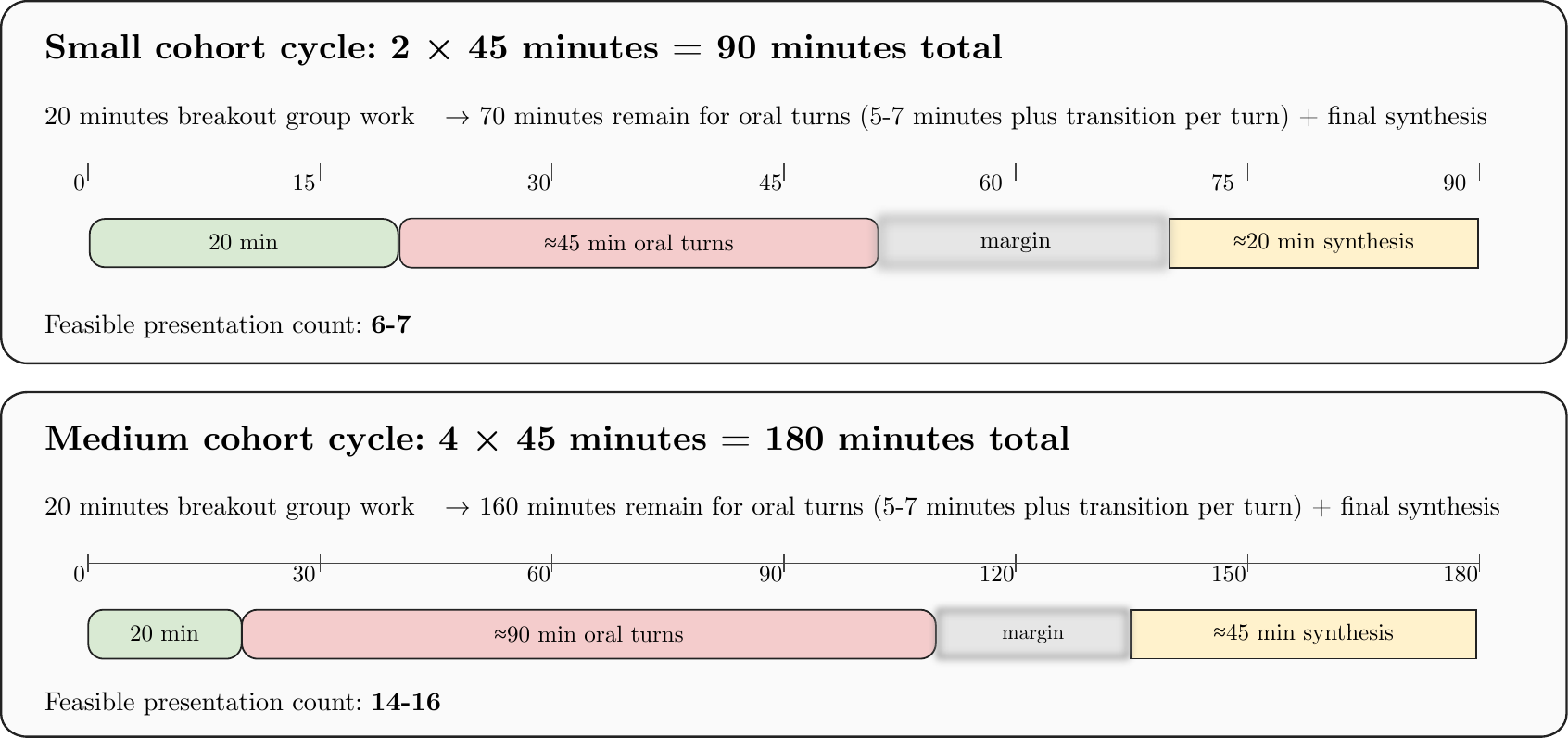}
\caption{\label{fig:orgc591a55}The figure visualizes the practical time arithmetic of the method for two common contact-time regimes.}
\end{figure*}

Table \ref{tab:org073247a} states the same logic compactly, and the 90 minute case is
illustrated in \ref{fig:orgc591a55}.

The model is therefore native for small and medium cohorts, up to about 60, provided
contact time grows accordingly. It does not scale to 300 by rhetoric alone. At that
size it requires smaller local groups and distributed human teaching capacity. This
does not make the method merely boutique. It means the general methodology is
human-scale even where institutions are large: large systems must be decomposed into
local groups if they want serious oral evidencing.

A practical minimum is at least three substantial primary observations per student
across the semester. In a 15-week semester with approximately four weekly
contact-hours, four observations becomes a theoretical upper target for cohorts
near 60. Additional late-semester sessions may be used to resolve borderline cases or
sparse visibility.

Where weekly visibility is too sparse, sampling-deficit recovery may come through
colloquia, oral checkpoints, or additional sampled observations in the same
cumulative system.

\begin{table*}[htbp]
\caption{\label{tab:org073247a}Human-time budget.}
\centering
\begin{tabularx}{\textwidth}{lYYY}
Aspect & Small cohort & Medium cohort & Large cohort\\
\hline
Contact group & about 30 students & about 60 students & 100+\\
Contact time & \(2 \times 45\) minutes & \(4 \times 45\) minutes & n/a\\
Group work & \(\approx\) 20 minutes & \(\approx\) 20 minutes & n/a at plenary scale\\
Presentation + synthesis & \(\approx\) 70 minutes & \(\approx\) 160 minutes & Not viable in one room\\
Supported presentations & \(\approx\) 6--8 gross task slots & up to 16 gross task slots & Requires smaller groups\\
Interpretation & Native case & Upper feasible range & Needs smaller local groups + teachers\\
\end{tabularx}
\end{table*}
\section{Fairness, evidence, and criterion-referenced grading}
\label{sec:org955f916}

The assessment logic begins from a simple distinction.

\begin{enumerate}
\item \emph{Primary evidence} comes from presentation plus follow-up questioning.
\item \emph{Supporting evidence} comes from substantial plenary intervention that clarifies,
repairs, extends, or redirects the mathematics.
\end{enumerate}

Table \ref{tab:orgfc9384d} formalizes the hierarchy.

\begin{table*}[htbp]
\caption{\label{tab:orgfc9384d}Evidence hierarchy.}
\centering
\begin{tabularx}{\textwidth}{lYl}
Evidence status & What it looks like & Counts how?\\
\hline
Primary & Designated presenter explains and survives follow-up & Direct grade-relevant evidence\\
Supporting & Useful question, clarification, correction, or tangent from a non-presenting student & Supplementary only\\
Promoted & Decisive correction / synthesis from a non-presenting student; survives brief follow-up & May become primary; justify in log\\
Does not count & Casual talk; repetition; dominance without mathematical value & Not grade-relevant\\
\end{tabularx}
\end{table*}

Supporting evidence can enrich judgment and help bridge sampling gaps, but it does
not normally replace primary evidence. It may, however, be promoted to primary status
when a student makes a mathematically illuminating contribution that is clearly tied
to the learning outcomes, survives brief follow-up questioning, and effectively
becomes the core explanation, correction, or synthesis in that moment. This is not
reward for extroversion; it is teacher-moderated recognition of a mathematically
decisive intervention. Whenever such a promotion occurs, it should be recorded with a
short justification note in the cumulative log.

The system is criterion-referenced rather than norm-referenced. Students are judged
against syllabus-level outcomes made visible in performance, not against cohort
rank. No curve-grading logic belongs here. The question is not whether a student
performed better than others, but what learning outcomes became visible, at what
quality, and through what kind of mathematical act.

The institutional weighting of oral work is local. Some frameworks allow more oral
time; some sharply limit it. This paper therefore does not fix a universal
percentage. What is constant is the validity carried by oral work in this method, and
also its constitutive role in learning. Oral work is the center of evidencing here
even where it is not 100\% of the grade. It serves as the strongest validity layer
over open-preparation artifacts and, where local rules permit, may also carry direct
grade weight. Local implementations may also include independent minority assessment
components without displacing the oral layer's central role.

In medium cohorts where observation density is imperfect, the cumulative oral record
can validate, challenge, and trigger further inquiry into the other evidence
available to the teacher. In that sense, the oral layer is not merely another
component; it is also the place where the trustworthiness of other evidence can be
tested.

The relevant comparison is not raw duration of observation. It is the number of
distinct opportunities in which understanding becomes visible under different local
conditions. A single three-hour written exam often produces one thick but noisy
terminal snapshot. Repeated shorter observations can reveal development, repair,
flexibility, and retained misunderstanding across time.

The pedagogical value for the non-presenting students also matters. Listening
students are not idle spectators. They are expected to compare the live explanation
with their own preparation, anticipate weaknesses, formulate questions, and
contribute supporting evidence where appropriate. The public explanation is therefore
simultaneously evidence for one student and a learning event for the rest of the
group.
\section{Five-grade oral rubric}
\label{sec:orgadb65cf}

The oral rubric should reflect visible outcome quality rather than whether a full
finished solution was presented. It should be read across three central axes:
mathematical correctness, formal expression, and responsive ownership under
questioning. A grade classification is given in Table \ref{tab:orgcef9b19}, while Table
\ref{tab:org7e16ebe} shows how one oral event can be recorded against outcomes.

\begin{table*}[htbp]
\caption{\label{tab:orgcef9b19}Five-grade oral evidence scale.}
\centering
\begin{tabularx}{\textwidth}{lYYY}
Grade & Correctness & Formal expression & Responsive ownership\\
\hline
5 Excellent & Correct; structurally sound & Disciplined; explicit & Adaptive; stable in follow-up\\
4 Very good & Correct; well-structured & Clear; mostly disciplined & Good response; minor limits\\
3 Satisfactory & Core idea correct but brittle & Partly controlled & Uneven or weak response\\
2 Barely passing & Minimal and fragile visibility & Imprecise; weakly organized & Heavy prompting needed\\
1 Not usable & No reliable outcome evidence & Missing or misleading & No meaningful survival\\
\end{tabularx}
\end{table*}

\begin{table*}[htbp]
\caption{\label{tab:org7e16ebe}Minimal outcome-linked micro-rubric.}
\centering
\begin{tabularx}{\textwidth}{lYY}
Dimension & Guiding question & Note type\\
\hline
Targeted outcome & Which learning outcome was visible? & Outcome code(s)\\
Correctness & Was the mathematics correct and non-misleading? & 1--5\\
Formal expression & Was the explanation coherent and rigorous? & 1--5\\
Ownership & Could the student adapt, justify, repair? & 1--5\\
Support & Was there a significant contribution off-turn? & Short comment; may supplement / promote\\
Next step & What should be checked next time? & Short comment\\
\end{tabularx}
\end{table*}

As a concrete example, take the grayscale dialogue from Listing
\ref{lst:grayscale-dialogue} and suppose it is judged against two outcomes:

\begin{itemize}
\item LO-G1: analyze and justify a linear transformation through basis images and
matrix-style reasoning;
\item LO-G2: interpret the structural effect of repeated application of an operator.
\end{itemize}

A student might receive the following paired judgment:

\begin{itemize}
\item LO-G1 = 2 (\emph{barely passing}), if the student recognizes that the operator is
determined by the images of the basis vectors, but can justify that only with
strong prompting and uncertain language;
\item LO-G2 = 3 (\emph{satisfactory}), if the student can explain that repeated application
changes nothing once the pixel is already gray, and can connect that fact to
projection-like or idempotent behavior in a basically correct way.
\end{itemize}

That is criterion-referenced judgment: not ``better than the class average,'' not
``presented a full solution,'' but ``these outcomes became visible at these quality
levels.''

A weak but barely passing performance might look like this:

\begin{quote}
The student correctly states that the grayscale procedure turns color into gray
and can compute output values for pure red, green, and blue, but initially treats
those computations as isolated examples rather than as images of basis vectors.
Under follow-up, the student eventually reaches the idea that knowing those three
images determines the whole map, but only after strong prompting and with hesitant
terminology. The student can also say that applying the operator twice does not
change an already gray pixel, but cannot yet explain this in a disciplined way as
projection-like or idempotent behavior. The core outcomes are minimally visible,
but the explanation remains fragile and only partly owned.
\end{quote}

This is the kind of performance that belongs around grade 2 for explanation and
perhaps grade 3 for a narrow structural outcome, depending on the exact outcomes
targeted.

Practical consistency in grading does not come from pretending oral judgment is
noise-free. It comes from anchored descriptors, immediate note-taking, and later
longitudinal review of accumulated evidence rather than one-shot impressionistic
judgment.
\section{Teacher load and cognitive load}
\label{sec:orgc62c982}

Teacher load is not a side issue. It is a design constraint. The method remains
viable only if both workload and cognitive load are intentionally managed. Table
\ref{tab:orgecc189f} states the distinction in compact form.

\begin{table*}[htbp]
\caption{\label{tab:orgecc189f}Teacher load in the model.}
\centering
\begin{tabularx}{\textwidth}{lYY}
Aspect & Workload & Cognitive load\\
\hline
Main source & Drafting; handouts; logging; tuning & Live orchestration; questioning; judging\\
Mitigation & AI support; bounded tasks; stable rubric model & Routine; stable criteria; small tasks; synthesis\\
\end{tabularx}
\end{table*}

The main mitigation is not heroic efficiency but structural consistency. Small tasks,
repeated routine, AI support, and a stable data-and-rubric model keep the system from
collapsing under its own good intentions. This also matters because teachers vary
widely in how they experience load. The same process that is manageable for one
instructor may become cognitively costly for another, including neurodiverse teachers
such as some teachers with ASD or ADHD. Easy adaptations therefore matter: more
explicit routine, a more stable note-taking template, smaller task bundles, or a more
strongly structured synthesis script may all preserve the method without changing its
core logic.
\section{Human-scale compactness}
\label{sec:orge55d138}

The method assumes that good mathematical teaching works with human-scale
compactness. A task, explanation, diagram, or code fragment should be small enough to
be held, traversed, and reconnected, but rich enough to alter the learner's web of
relations. Understanding grows when an idea is seen in many contexts, tested on
examples, and tied to several metaphors until it changes how neighboring topics are
understood. The educational point is not maximal symbolic density. It is compactness
that can be unfolded. This idea is close to Greg Egan's depiction of mathematical
understanding as a dense web of mutually illuminating relations among symbols and
ideas \citep{Egan1997Diaspora}, and it sits naturally beside the anti-helicopter
warning in current AI discussion: the point is not to be dropped at the destination,
but to traverse enough structure that one's thinking changes \citep{Wong2026TaoAI}.

This is one reason the method resists pure push-button tool use. A point-and-click
workflow can produce a result quickly while leaving the relational web untouched. The
aim is not to abolish tools, but to ensure that students learn functionalities and
structures they can later compose freely, rather than becoming dependent on one
narrow mediation path.
\section{Discussion and short closing remarks}
\label{sec:org8e3e266}

ABC Learning Design is relevant here as one example of a neighboring design
grammar. It offers a storyboard-oriented way of thinking about course design through
recurring types of learning activity. That matters in the present context because the
current paper also proposes a recurring grammar, though at a different scale and with
different emphases. The weekly cycle developed here can likewise be read through
activity types such as acquisition, investigation, discussion, practice,
collaboration, and production, but with oral validation and instructor synthesis
carrying unusually central weight \citep{ABC2020}. In that sense, ABC is useful here
less as a theory than as a neighboring grammar for variation.

A broader mathematics-community discussion about AI is also worth noting. The AMS has
recently begun organizing advisory and white-paper work on AI and mathematics, which
reflects the same general concern visible here: standards, norms, and workflows are
being renegotiated quickly, and mathematics education cannot simply wait for the dust
to settle \citep{AMS2026AI}.

The manuscript itself has also been structured to reflect its own methodology in a
restrained way. It moves from introduction, to a short concrete example, to later
conceptual exposition. That is not meant as a gimmick. It is meant to let the method
show some of its own plausibility in the very form of the article.

A related question lies just beyond the present paper: what a similar system might
look like for online part-time students. That is both a limit and an opportunity. The
present article develops the method for in-person human-scale teaching, but the
underlying logic may extend further if one preserves small-group oral ownership,
contingent questioning, and cumulative recording while adapting scheduling, breakout
structure, and asynchronous preparation to online constraints.

This remains a methodology paper rather than an effectiveness study. The next
empirical questions are therefore local and manageable: how many primary observations
are actually achieved under different contact structures; which task families best
produce owned explanation rather than memorized script; how much the two Oracles
reduce teacher load without over-scaffolding students; how students experience the
transition from bootstrapping to accountability; and which forms of realism most
consistently generate explanations that survive follow-up.
\section{Conclusion}
\label{sec:org3a69c1b}

Under current AI conditions, the central assessment question is no longer whether
students can submit polished artifacts. It is where trustworthy evidence of
understanding should now be sought. This paper proposes one answer. Use short
coherent problems. Let students prepare openly. Require human explanation. Follow it
with contingent questioning. End with synthesis. Record the evidence cumulatively.

That answer is neither anti-tool nor anti-lecture. It is an RME-guided,
question-first, orally validated methodology designed for AI-compatible higher
education. Its immediate purpose is practical: to make understanding visible without
pretending that the technological environment has not changed. It is also one
alternative to the current drift toward proctored digital exams: instead of
escalating surveillance, it relocates trust toward human explanation, contingent
questioning, and cumulative judgment.
\section{Declaration}
\label{sec:org591b046}

The author reports no competing or conflicts of interest regarding this
manuscript. This work was supported by the Faculty of Informatics of Juraj Dobrila
University of Pula.

\bibliographystyle{authoryear}
\bibliography{references1}

\begin{thebibliography}{23}
\providecommand{\natexlab}[1]{#1}
\providecommand{\url}[1]{\texttt{#1}}
\expandafter\ifx\csname urlstyle\endcsname\relax
  \providecommand{\doi}[1]{doi: #1}\else
  \providecommand{\doi}{doi: \begingroup \urlstyle{rm}\Url}\fi

\bibitem[{ABC Learning Design}(2020)]{ABC2020}
{ABC Learning Design}.
\newblock Abc learning design, 2020.
\newblock URL \url{https://abc-ld.org/}.
\newblock Storyboard-based learning design workshop and toolkit.

\bibitem[{American Mathematical Society}(2026)]{AMS2026AI}
{American Mathematical Society}.
\newblock Artificial intelligence and mathematics, 2026.
\newblock URL \url{https://www.ams.org/about-us/ai}.
\newblock AMS page aggregating advisory-group work, white papers, videos, and
  articles on AI and mathematics.

\bibitem[Baartman et~al.(2022)Baartman, van Schilt-Mol, and van~der
  Vleuten]{BaartmanVanSchiltMolVanderVleuten2022}
Liesbeth Baartman, Tamara van Schilt-Mol, and Cees van~der Vleuten.
\newblock Programmatic assessment design choices in nine programs in higher
  education.
\newblock \emph{Frontiers in Education}, 7:\penalty0 931980, 2022.
\newblock \doi{10.3389/feduc.2022.931980}.

\bibitem[Biggs(1996)]{Biggs1996Constructive}
John Biggs.
\newblock Enhancing teaching through constructive alignment.
\newblock \emph{Higher Education}, 32\penalty0 (3):\penalty0 347--364, 1996.
\newblock \doi{10.1007/BF00138871}.

\bibitem[Egan(1997)]{Egan1997Diaspora}
Greg Egan.
\newblock \emph{Diaspora}.
\newblock Orion/Millennium, London, 1997.
\newblock ISBN 1-85798-438-2.

\bibitem[Freudenthal(1991)]{Freudenthal1991Revisiting}
Hans Freudenthal.
\newblock \emph{Revisiting Mathematics Education: China Lectures}.
\newblock Kluwer Academic Publishers, Dordrecht, 1991.

\bibitem[{Google for Education}(2026)]{GoogleNotebookLM2026}
{Google for Education}.
\newblock Understand anything with notebooklm, 2026.
\newblock URL \url{https://edu.google.com/ai-notebooklm/}.
\newblock Official NotebookLM for education page describing grounded use over
  uploaded sources.

\bibitem[Gravemeijer and Doorman(1999)]{GravemeijerDoorman1999}
Koeno Gravemeijer and Michiel Doorman.
\newblock Context problems in realistic mathematics education: A calculus
  course as an example.
\newblock \emph{Educational Studies in Mathematics}, 39\penalty0
  (1--3):\penalty0 111--129, 1999.
\newblock \doi{10.1023/A:1003749919816}.

\bibitem[Hmelo-Silver(2004)]{HmeloSilver2004PBL}
Cindy~E. Hmelo-Silver.
\newblock Problem-based learning: What and how do students learn?
\newblock \emph{Educational Psychology Review}, 16\penalty0 (3):\penalty0
  235--266, 2004.
\newblock \doi{10.1023/B:EDPR.0000034022.16470.f3}.

\bibitem[Johnson(1993)]{Johnson1993Timetabling}
Douglas Johnson.
\newblock A database approach to course timetabling.
\newblock \emph{INFORMS Journal on Computing}, 5\penalty0 (3):\penalty0
  302--310, 1993.

\bibitem[Kurnik(2008)]{Kurnik2008Znanstvenost}
Zdravko Kurnik.
\newblock Znanstvenost u nastavi matematike.
\newblock \emph{Metodika}, 9\penalty0 (17):\penalty0 318--327, 2008.
\newblock URL \url{https://hrcak.srce.hr/34802}.

\bibitem[Kurnik(2009)]{Kurnik2009Problemska}
Zdravko Kurnik.
\newblock Problemska nastava.
\newblock \emph{Mi{\v S} -- Matematika i {\v S}kola}, 2009.
\newblock URL \url{https://mis.element.hr/clanak/problemska-nastava/}.

\bibitem[Laursen and Rasmussen(2019)]{LaursenRasmussen2019IBL}
Sandra~L. Laursen and Chris Rasmussen.
\newblock I on the prize: Inquiry approaches in undergraduate mathematics.
\newblock \emph{International Journal of Research in Undergraduate Mathematics
  Education}, 5\penalty0 (1):\penalty0 129--146, 2019.
\newblock \doi{10.1007/s40753-019-00085-6}.

\bibitem[Li et~al.(2025)]{LiEtAl2025RAGEducation}
Z.~Li et~al.
\newblock Retrieval-augmented generation for educational application: A survey.
\newblock \emph{Smart Learning Environments}, 2025.
\newblock Survey discussing RAG in educational settings and ongoing challenges
  such as hallucination, outdated knowledge, and reliability limits.

\bibitem[Lo et~al.(2017)Lo, Hew, and Chen]{LoHewChen2017MathFlipped}
Chung~Kwan Lo, Khe~Foon Hew, and Gaowei Chen.
\newblock Toward a set of design principles for mathematics flipped classrooms:
  A synthesis of research in mathematics education.
\newblock \emph{Educational Research Review}, 22:\penalty0 50--73, 2017.
\newblock \doi{10.1016/j.edurev.2017.08.002}.

\bibitem[Mac~Lane(1998)]{MacLane1998Categories}
Saunders Mac~Lane.
\newblock \emph{Categories for the Working Mathematician}, volume~5 of
  \emph{Graduate Texts in Mathematics}.
\newblock Springer, New York, 2 edition, 1998.
\newblock ISBN 978-0-387-98403-2.

\bibitem[M{\"u}ller and
  Bart{\'a}k(2002)]{MullerBartak2002InteractiveTimetabling}
Tom{\'a}{\v{s}} M{\"u}ller and Roman Bart{\'a}k.
\newblock Interactive timetabling: Concepts, techniques, and practical results.
\newblock \emph{Practice and Theory of Automated Timetabling}, 2002.

\bibitem[Pierce(1991)]{Pierce1991BasicCategoryTheory}
Benjamin~C. Pierce.
\newblock \emph{Basic Category Theory for Computer Scientists}.
\newblock Foundations of Computing. MIT Press, Cambridge, MA, 1991.
\newblock ISBN 978-0-262-66071-6.

\bibitem[Strelan et~al.(2020)Strelan, Osborn, and
  Palmer]{Strelan2020FlippedMeta}
Peter Strelan, Amanda Osborn, and Edward Palmer.
\newblock The flipped classroom: A meta-analysis of effects on student
  performance across disciplines and education levels.
\newblock \emph{Educational Research Review}, 30:\penalty0 100314, 2020.
\newblock \doi{10.1016/j.edurev.2020.100314}.

\bibitem[van Alten et~al.(2019)van Alten, Phielix, Janssen, and
  Kester]{VanAlten2019FlippedMeta}
David C.~D. van Alten, Chris Phielix, Jeroen Janssen, and Liesbeth Kester.
\newblock Effects of flipping the classroom on learning outcomes and
  satisfaction: A meta-analysis.
\newblock \emph{Educational Research Review}, 28:\penalty0 100281, 2019.
\newblock \doi{10.1016/j.edurev.2019.05.003}.

\bibitem[van~den Heuvel-Panhuizen and
  Drijvers(2014)]{HeuvelPanhuizenDrijvers2014RME}
Marja van~den Heuvel-Panhuizen and Paul Drijvers.
\newblock Realistic mathematics education.
\newblock In Stephen Lerman, editor, \emph{Encyclopedia of Mathematics
  Education}, pages 521--525. Springer, Dordrecht, 2014.
\newblock \doi{10.1007/978-94-007-4978-8_170}.

\bibitem[van~der Vleuten et~al.(2012)van~der Vleuten, Schuwirth, Driessen,
  Dijkstra, Tigelaar, Baartman, and van
  Tartwijk]{VanderVleutenEtAl2012Programmatic}
Cees P.~M. van~der Vleuten, Lambert W.~T. Schuwirth, Erik~W. Driessen, Janke
  Dijkstra, Dineke Tigelaar, Liesbeth K.~J. Baartman, and Jan van Tartwijk.
\newblock A model for programmatic assessment fit for purpose.
\newblock \emph{Medical Teacher}, 34\penalty0 (3):\penalty0 205--214, 2012.
\newblock \doi{10.3109/0142159X.2012.652239}.

\bibitem[Wong(2026)]{Wong2026TaoAI}
Matteo Wong.
\newblock The edge of mathematics.
\newblock \emph{The Atlantic}, February 2026.
\newblock URL
  \url{https://www.theatlantic.com/technology/2026/02/ai-math-terrance-tao/686107/}.

\end{thebibliography}
\end{document}